\documentclass[]{interact}
\usepackage{graphicx}
 \usepackage{cite}
\usepackage{amsmath}
\usepackage{float}
\usepackage{tabularx}
\usepackage{array}
\usepackage[font=small,labelfont=bf]{caption}
\usepackage{graphicx}
\usepackage{xcolor}
\usepackage{rotating,booktabs,multirow}
\usepackage{varwidth}
\usepackage{epstopdf}
\usepackage[caption=false]{subfig}

\usepackage[numbers,sort&compress]{natbib}
\bibpunct[, ]{[}{]}{,}{n}{,}{,}

\theoremstyle{plain}
\newtheorem{theorem}{Theorem}[section]

\newtheorem{proposition}[theorem]{Proposition}
\newcommand{\argmin}{\operatornamewithlimits{argmin}}
\DeclareMathOperator{\conv}{conv}
\DeclareMathOperator{\diag}{diag}
\DeclareMathOperator{\Diag}{Diag}
\DeclareMathOperator{\Int}{int}
\theoremstyle{definition}

\newtheorem{example}[theorem]{Example}

\theoremstyle{remark}

\begin{document}


\title{New bounds for nonconvex quadratically constrained quadratic programming}

\author{
\name{Moslem  Zamani\textsuperscript{a}\textsuperscript{b}\textsuperscript{c}\thanks{ M. Zamani. Email: zamani.moslem@tdt.edu.vn}}
\affil{\textsuperscript{a} Parametric MultiObjective Optimization Research Group, Ton Duc Thang University, Ho Chi Minh City, Vietnam; 
\textsuperscript{b} Faculty of Mathematics and Statistics, Ton Duc Thang University, Ho Chi Minh City, Vietnam; 
\textsuperscript{c} School of Mathematics, Statistics and Computer Science, College of Science, University of Tehran, Enghelab Avenue, Tehran, Iran; }
}

\maketitle

\begin{abstract}
In this paper, we study some bounds for nonconvex quadratically constrained quadratic programs. We propose two types of bounds for quadratically constrained quadratic programs, quadratic and cubic bounds. For quadratic bounds, we use affine functions as Lagrange multipliers. We demonstrate that  most semi-definite relaxations can be obtained as the dual of  a quadratic bound. In addition, we study bounds obtained by changing the ground set. For cubic bounds, in addition to affine multipliers we employ quadratic functions. We provide a comparison between the proposed cubic bound and typical bounds for standard quadratic programs. Moreover, we report comparison results of a quadratic and a cubic bound  for  some
 non-convex quadratically constrained quadratic programs.
\end{abstract}

\begin{keywords}
Quadratically constrained quadratic programming; Semidefinite relaxation; Reformulation-linearization technique 
\end{keywords}
\section{Introduction}
\noindent
 We consider the following  quadratically constrained quadratic programming, QCQP,
\begin{equation}\tag{QCQP}\label{P}
\begin{array}{ll} 
 \ & \min \  x^TQ_0x+2c_0^Tx \\ 
 &  s.t. \  x^TQ_ix+2c_i^Tx\leq b_i, \ \ i=1, ..., m \\
   & \; \; \; \; \; \; \; \; \; Ax=d \\  
  & \; \; \; \; \; \; \; \; \;   l\leq x\leq u,
\end{array}
\end{equation}
where $x\in\mathbb{R}^n$ is the vector of decision variables, $Q_i$ ($i=0, 1, ..., m$) are $n\times n$  real symmetric matrices, $A$ is  a $p\times n$  real matrix, $c_i$ ($i=0,1, ..., m$) and $d$ are vectors in $\mathbb{R}^n$ and $\mathbb{R}^p$, respectively,  and $b_i$ ($i=1, ..., m$) are real scalars. We assume that $-\infty < l_i\leq u_i < \infty $ for $i=1, ..., n$. Without loss of generality, we may assume that $l=0$ and $u=e$, where $e$ represents vector of ones in $\mathbb{R}^n$. We remark that general QCQPs with bounded feasible set can be formulated as \eqref{P}.\\ 
\indent
QCQP is a fundamental problem in optimization theory and practice. QCQPs  arise in many applications including economic equilibria, facility location and circle packing problems \cite{Bao, pack, fac}. Furthermore, most combinatorial optimization problems including max-cut problem and clique problem can be casted as QCQPs \cite{Goemans, Cli}. In addition to the aforementioned problem, Madani et al. showed that any polynomial optimization problem can be casted as a QCQP  \cite{Madani}. When the matrices $Q_i$ ($i=0, 1, ..., m$) are positive semi-definite \eqref{P} will be a convex optimization problem, and it is polynomially solvable. Nevertheless, as QCQPs include a wide range of  NP-hard optimization problems, QCQP is NP-hard \cite{Va}.\\
\indent
A typical class of optimization  methods for handling QCQPs is branch-and-bound method. In this approach, the main problem is divided to some subproblems, which are called nodes. At each node, a lower bound is computed by a relaxation or a bound. In general, the generated lower bound determines  a node will be fathomed or branched. The effectiveness of a branch-and-bound method rests mainly on the tightness of generated lower bounds and their computational time. Most relaxation approaches and bounds for QCQPs are mainly based on  the reformulation-linearization technique (RLT), convex relaxations and semidefinite relaxations \cite{Tuy, Sherali, Locatelli}. The most effective relaxation method is semi-definite relaxation (SDR) \cite{Bao, Koj}.  Due to the efficiency of this approach, many SDRs have been proposed; See \cite{Bao, Zheng} for review and comparisons. \\
\indent
Recently, Zamani has proposed a new dual for linearly constrained quadratic programming \cite{Za}. He considers affine functions as Lagrange multipliers. In this paper, similar to his method, we present two types of bounds for QCQPs, quadratic and cubic bounds. For quadratic bounds, we employ affine functions as Lagrange multipliers. We illustrate that most SDRs can be interpreted as the dual of a quadratic bound. In addition, we introduce  some bounds which are obtained by changing the ground set. \\
\indent 
 For cubic bounds, we apply quadratic functions as Lagrange multipliers. We give some conditions under which the proposed bound is exact. We  demonstrate that the cubic bound is equivalent the bound obtained by Parrilo hierarchy for standard quadratic programs. \\
\indent
 The paper is organized as follows. After reviewing our notations, in Section 2 we introduce the quadratic bounds. Section 3 is devoted to cubic bounds.  In Section 4, we illustrate the effectiveness of a quadratic and a cubic bound by presenting its numerical performance on some  QCQPs.
 \subsection{Notation}
The following notation is used throughout the paper. The $n$-dimensional Euclidean
space is denoted by $\mathbb{R}^n$. Let $A_i$ stand for the $i^{th}$ row of matrix $A$. Vectors are considered to be column vectors and $T$ denotes transposition operation. We employ $e_i$ to represent  the $i^{th}$ unit vector, and  vector $e$ stands for vector of ones. 
For symmetric matrices $A$ and $B$, we use notation $A \succeq B$ to denote $A-B$ is positive semidefinite. The inner
product $A$ and $B$ is defined and denoted as $A \bullet B=trace(AB)$. A symmetric $n\times n$ matrix $Q$ is called copositive if bilinear form $x^TQx$ is non-negative on non-negative orthant. For $x\in\mathbb{R}^n$, we denote by $\diag(x)$ the diagonal matrix with $\diag(x)_{ii}=x_i$. Moreover,  for $n\times n$ matrix $Q$, 
$\Diag(Q)$ denotes a column vector with $\Diag(Q)_i=Q_{ii}$.\\
\indent
For a set $\mathcal{X}\subseteq\mathbb{R}^n$, we use the notations $\Int(\mathcal{X})$ and
$\conv(\mathcal{X})$ for the interior and the convex hull of  $\mathcal{X}$, respectively. $\mathbb{R}_{+}^n$ denotes non-negative orthant. We use $B$ to represent box $[0, 1]^n$. The dual cone of  $K$ is denoted by $K^*=\{y: y^Tx\geq 0, \ \forall x\in K\}$.\\
\indent
 We use $\mathcal{A}(\mathbb{R}^n)$ and $\mathcal{Q}(\mathbb{R}^n)$ to represent affine and quadratic functions on $\mathbb{R}^n$. We denote non-negative affine and quadratic functions on  $\mathcal{X}\subseteq \mathbb{R}^n$ by $\mathcal{A}_+(\mathcal{X})$ and $\mathcal{Q}_+(\mathcal{X})$, respectively, that is 
 $\mathcal{A}_+(\mathcal{X})=\{\alpha\in\mathcal{A}(\mathbb{R}^n): \alpha(x)\geq 0 \ \forall x\in \mathcal{X}\}$ and $\mathcal{Q}_+(\mathcal{X})=\{q\in\mathcal{Q}(\mathbb{R}^n): q(x)\geq 0 \ \forall x\in \mathcal{X}\}$.
  For a quadratic function $q(x) = x^T Qx +2c^T x + b$, we
denote the matrix representation of  $q$ by $\mathcal{M}(q)=\begin{pmatrix} Q & c \\ c^T & b \end{pmatrix}$.

\section{Quadratic bounds}
\label{sec:1}
\noindent
 In this section, we  propose some quadratic bounds for QCQP. Let $\mathcal{X}=\{x\in \mathbb{R}^n: x\in B, Ax=d\}$ and 
 $F=\{x\in \mathbb{R}^n: x\in \mathcal{X}, x^TQ_ix+2c_i^Tx\leq b_i,  i=1, ..., m \}$. Because $\mathcal{A}_+(\mathcal{X})$ is a polyhedral set, it follows that the representation of $\mathcal{A}_+(\mathcal{X})$ in $\mathbb{R}^{n+1}$ is a polyhedral cone \cite{Mang}. We propose the following problem as a new quadratic bound for  \eqref{P}, 
\begin{align} \label{L2}
\nonumber \ & \max \ \ell \\
\nonumber  \   &  s.t.\   x^TQ_0x+2c^T_0x-\ell+\sum_{i=1}^{m} \lambda_i(x^TQ_ix+2c^T_ix-b_i)+ \sum_{i=1}^{p} \alpha_i(x)(A_ix-d_i)-\\
  & \; \; \; \; \; \; \; \; \; \; \; \;         \sum_{i=1}^{n} \beta_i(x)x_i+\sum_{i=1}^{n} \gamma_i(x)(x_i-1)\in \mathcal{Q}_+(\mathbb{R}^{n}),\\
\nonumber & \; \; \; \; \; \; \; \; \; \; \; \;     \ \alpha_i\in \mathcal{A}(\mathbb{R}^n), \ i=1, ..., p\\
\nonumber & \; \; \; \; \; \; \; \; \; \; \; \;   \lambda_i \geq 0,\ \beta_i, \gamma_i \in \mathcal{A}_+(\mathcal{X}), \ i=1, ..., n
\end{align}
 The above problem can be regarded as a Lagrangian dual for  \eqref{P} for which the dual variables corresponding to linear constraints are replaced by affine functions. We remark that, due to the non-homogeneous Farkas' Lemma,  $\alpha(x)=f^Tx+g$ belong to $\mathcal{A}_+(\mathcal{X})$ if and only if there exist $\lambda\in\mathbb{R}^p$ and  $\mu\in\mathbb{R}_{+}^n$  with
 $f\geq A^T\lambda-\mu$ and $g\geq d^T\lambda+e^T\mu$.
  Note that the quadratic function $q(x)=x^TQx+2c^Tc+b$ is non-negative on $\mathbb{R}^{n}$ if and only if matrix
$\mathcal{M}(q)$ is positive semi-definite, and accordingly
  problem \eqref{L2} can be formulated as a semi-definite program, which has $O(n^2)$ variables. \\
  \indent
One crucial question regarding this bound is well-definedness.  In next proposition, we prove that problem \eqref{L2} is feasible and generates a finite lower bound for  \eqref{P}.
\begin{proposition}\label{p1}
Let \eqref{P} have a feasible point. Then problem \eqref{L2}  gives a finite lower bound.
\begin{proof} 
Similar to the proof of Proposition 2 in \cite{Za}, it is shown that there exist $\gamma_i \in \mathcal{A}_+(\mathcal{X})$ for  $i=1, ..., n$ such that  $x^TQ_0x+2c^T_0x+\sum_{i=1}^{n} \gamma_i(x)(x_i-1)$ is strictly convex. So for suitable choice of $\ell$, we have  $x^TQ_0x+2c^T_0x+\sum_{i=1}^{k} \gamma_i(x)(x_i-1)-\ell\in \mathcal{Q}_+(\mathbb{R}^{n})$, which shows the feasibility of \eqref{L2}. The first constraint of  \eqref{L2} implies that the optimal value of \eqref{L2} is a lower bound for \eqref{P}.
\end{proof} 
\end{proposition}
The proof of above theorem reveals that problem \eqref{L2} is feasible for each quadratic function. In fact, the problem is strongly feasible.  The conic optimization problem is called strongly feasible if it is feasible and remains feasible for all sufficiently small perturbations of right side of linear constraints \cite{Ren}. 
 As problem  \eqref{L2}  is convex, it is natural to ask about the dual thereof. The  dual of  \eqref{L2} can be written as 
 \begin{align}\label{D2}
\nonumber  \min \ & Q_0\bullet X+2c_0^Tx \\
\nonumber      s.t.\ &  Q_i\bullet X+2c_i^Tx\leq b_i, \ \ \;\;\;\;\;\;\;\;\;\  i=1, ..., m \\
\nonumber      &  XA_i^T= d_ix, \ \ \;\;\;\;\;\;\;\;\;\;\;\;\; \;\;\;\;\;\;\;\ i=1, ..., p \\
\nonumber    & Ax=d,  \\
&    X\geq 0, \\
\nonumber  &   ex^T-X\geq 0,\\
\nonumber  &   X-ex^T-xe^T+ee^T\geq 0,\\
\nonumber &   X-xx^T\succeq 0.
\end{align}
We refer the reader for the details of computation of \eqref{D2} to \cite{Za}. In problem \eqref{D2}, we did not write redundant constraint $x\in B$ \cite{Sher}. Problem \eqref{D2} is a well-known Shor relaxation with partial first-level RLT \cite{Ans, Bao}. Anstreicher proposed SDR \eqref{D2} as a combination of  RLT and Shor relaxation. He showed that SDR \eqref{D2} can generate bounds tighter than either technique. Bao et al. established that SDR \eqref{D2} and  doubly non-negative relaxation provide the same bound. The doubly non-negative relaxation is similar to problem \eqref{D2}, but the constraint $XA_i^T= d_ix$ is replaced by $X\bullet A_i^TA_i= d_i^2$, $i=1, ..., m $. 
Note that as \eqref{L2} is strongly feasible, strong duality holds \cite{Ren}, and consequently problems \eqref{L2} and \eqref{D2} are equivalent.\\
\indent
 Since the ground set of \eqref{P} is not $\mathbb{R}^n$, the bound may be improved if we replace $\mathcal{Q}_+(\mathbb{R}^{n})$ by other sets. Bomze took advantage of this idea and has proposed some results about global optimality conditions for QCQPs \cite{Bom}.
As the feasible set of \eqref{P} is subset of positive orthant, one replacement for  $\mathcal{Q}_+(\mathbb{R}^{n})$ can be quadratic functions with non-negative coefficients. In this case, we have the following bound
 \begin{align}\label{L3}
\nonumber  \max \ & \ell \\
\nonumber   s.t.\    &  x^TQ_0x+2c^T_0x-\ell+\sum_{i=1}^{m} \lambda_i(x^TQ_ix+2c^T_ix-b_i)+ \\
\nonumber  & \; \;  \; \;        \sum_{i=1}^{p} \alpha_i(x)(A_ix-d_i)-\sum_{i=1}^{n} \beta_i(x)x_i+\sum_{i=1}^{n} \gamma_i(x)(x_i-1)\in
 \mathcal{Q}^N(\mathbb{R}^{n}),\\
& \lambda_i\geq 0, \ \alpha_i\in \mathcal{A}(\mathbb{R}^n), \ i=1, ..., m,\\
\nonumber &   \beta_i, \gamma_i \in \mathcal{A}_+(\mathcal{X}), \ i=1, ..., p.
\end{align}
where $\mathcal{Q}^N(\mathbb{R}^{n})$ denotes quadratic functions with non-negative coefficients. The above problem can be formulated as  a linear program. It can be shown that  \eqref{L3} is the dual of a linear RLT \cite{Sher}.\\
\indent
Another interesting substitute for  $\mathcal{Q}_+(\mathbb{R}^{n})$ is non-negative quadratic functions on $B$.  In this case, the following  program provides a bound
 \begin{align}\label{L4}
\nonumber  \max \ &\ell \\
\nonumber    s.t.\  & x^TQ_0x+2c^T_0x-\ell+\sum_{i=1}^{m} \lambda_i(x^TQ_ix+2c^T_ix-b_i)+   \sum_{i=1}^{p} \alpha_i(x)(A_ix-d_i)\in  \mathcal{Q}_+(B),\\
&    \lambda\geq 0, \ \alpha_i\in \mathcal{A}(\mathbb{R}^n) \  \ \ \ i=1, ..., m
\end{align}
Since for each $q\in \mathcal{Q}(\mathbb{R}^{n})$, there exists $\ell$ with $q-\ell\in\mathcal{Q}_+(B)$, problem \eqref{L4} is always feasible. Needless to say, bound \eqref{L4} dominates all the above-mentioned bounds. Nevertheless, this bound is not necessarily exact for general QCQPs. A bound or relaxation is said to be exact if it provides a bound equal to the optimal value of main problem.  Next theorem provides some sufficient conditions for exactness. 
\begin{theorem}\label{T1}
Bound \eqref{L4} is exact if 
 \begin{align*}
x^TQ_ix+2c^T_ix-b_i\leq (or \geq) 0, \  \  \forall x\in X, i=1, ..., m
\end{align*}
\begin{proof} 
First we prove the case that there does not exist any quadratic constraint. By Lemma 4 in \cite{Eish}, 
\begin{align*}
\mathcal{Q}_+(B)=\{\begin{pmatrix} x^T & 1 \end{pmatrix}Q \begin{pmatrix} x \\ 1 \end{pmatrix}: Q\in K_B^*\}\\
\mathcal{Q}_+(\mathcal{X})=\{\begin{pmatrix} x^T & 1 \end{pmatrix}Q \begin{pmatrix} x \\ 1 \end{pmatrix}: Q\in K_X^*\}
\end{align*}
where 
\begin{align*}
& K_B=\conv\{zz^T: z\in\mathbb{R}^{n+1}_+ , z_i\leq z_{n+1} \  i=1, ...,n\}\\
& K_\mathcal{X}=\conv\{zz^T: z\in\mathbb{R}^{n+1}_+ ,
 \begin{pmatrix} A & -b\end{pmatrix} z=0, \  z_i\leq z_{n+1} \  i=1, ...,n\}. 
 \end{align*}
 Remark that $ \mathcal{M}(\mathcal{Q}_+(B))^*=K_B$ and $ \mathcal{M}(\mathcal{Q}_+(\mathcal{X}))^*=K_\mathcal{X}$.
 We show that
 $\mathcal{Q}_+(\mathcal{X})=\mathcal{Q}_+(B)+\{ \sum_{i=1}^{p} \alpha_i(x)(A_ix-d_i): \alpha_i\in \mathcal{A}(\mathbb{R}^n)\}$. The inclusion $\supseteq$ is trivial. We prove the inclusion $\subseteq$ by contradiction. Let 
 $q(x)=x^TQx+2c^Tx+c_0\in \mathcal{Q}_+(\mathcal{X})$ 
 while $q\notin \mathcal{Q}_+(B)+\{ \sum_{i=1}^{p} \alpha_i(x)(A_ix-d_i): \alpha_i\in \mathcal{A}(\mathbb{R}^n)\}$.
  By separation theorem, there is $O\in \mathcal{M}(\mathcal{Q}_+(B))^*+\mathcal{M}(\{ \sum_{i=1}^{p} \alpha_i(x)(A_ix-d_i): \alpha_i\in \mathcal{A}(\mathbb{R}^n)\})^*$ with $\mathcal{M}(q)\bullet O=-1$. As $\mathcal{M}(\mathcal{Q}_+(B))^*=K_B$, we have $O=\sum_{k=1}^l z^k(z^k)^T$, where 
  $z^k\in\{\mathbb{R}^{n+1}_+ , z_i\leq z_{n+1} \  i=1, ...,n\}$ ($k=1, ..., l$). As $(a_ix-b)(a_ix-b)$ and $-(a_ix-b)(a_ix-b)$ are members of $\{ \sum_{i=1}^{p} \alpha_i(x)(A_ix-d_i): \alpha_i\in \mathcal{A}(\mathbb{R}^n)\}$ for $i=1, ..., p$, we have $O\in K_\mathcal{X}$. This implies that $O\in \mathcal{M}(\mathcal{Q}_+(\mathcal{X}))^*$. Thus, we have  $\mathcal{M}(q)\bullet O\geq 0$ which contradicts the assumption $\mathcal{M}(q)\bullet O=-1$.\\
  Now we consider the case that quadratic constraints exist. If $x^TQ_ix+2c^T_ix\leq b_i$ for each $x\in \mathcal{X}$ and $i=1, ..., m$,  the quadratic constraints are redundant  and theorem follows from the first part.\\
   For the case that $x^TQ_ix+2c^T_ix\geq b_i$ for each $x\in \mathcal{X}$ and $i=1, ..., m$, we establish that the dual cones corresponding to 
  $\mathcal{M}(\mathcal{Q}_+(F))$ and $\mathcal{M}(\mathcal{Q}_+(A)+\{ \sum_{i=1}^{m} \lambda_i(x^TQ_ix+2c^T_ix-b_i): \lambda\leq 0\})$ are the same. The inclusion 
  $\mathcal{M}(\mathcal{Q}_+(F))^*\subseteq \mathcal{M}(\mathcal{Q}_+(A)+\{ \sum_{i=1}^{m} \lambda_i(x^TQ_ix+2c^T_ix-b_i): \lambda\leq 0\})^*$
   is immediate. We show the reverse inclusion.  Suppose that  
 $O\in \mathcal{M}(\mathcal{Q}_+(A)+\{ \sum_{i=1}^{m} \lambda_i(x^TQ_ix+2c^T_ix-b_i): \lambda\leq 0\})^*$.
  By the representation of $\mathcal{Q}_+(A)$, we have $O=\sum_{k=1}^l z^k(z^k)^T$, where 
  $z^k\in\{\mathbb{R}^{n+1}_+ , \begin{pmatrix} A & -d \end{pmatrix}z=0, z_i\leq z_{n+1} \  i=1, ...,n\}$ ($k=1, ..., l$). By the assumption $\begin{pmatrix} Q_i & c_i\\ c_i^T & b_i \end{pmatrix}\bullet z_kz_k^T\geq 0$ ($i=1, ..., m, k=1, ..., l$), we have $z_i^TOz_i=0$. Suppose that 
\begin{align*}
K_F=\conv\{zz^T: z  \in&\mathbb{R}^{n+1}_+ ,  \begin{pmatrix} A & -d \end{pmatrix}z=0, z_i\leq z_{n+1} \  i=1, ...,n, z^T\mathcal{M}(q_j) z\leq 0,\\
&  j=1, ..., m\}.
\end{align*}
 As 
$\mathcal{Q}_+(F)=\{\begin{pmatrix} x^T & 1 \end{pmatrix}Q \begin{pmatrix} x \\ 1 \end{pmatrix}: Q\in K_F^*\}$,  we have $O\in \mathcal{M}(\mathcal{Q}_+(F))^*$, which completes the proof.
\end{proof} 
\end{theorem}
It is worth mentioning that the Theorem \ref{T1} can be proved by using strong duality for conic programs and Proposition 6 in \cite{Bao}, but here we present a new proof. 
Next proposition states that bound \eqref{L4} is exact for linearly constrained quadratic programs with binary variables.
\begin{proposition}
Bound \eqref{L4} is exact for linearly constrained quadratic programs with  binary variables.
\begin{proof} 
Consider the problem
\begin{equation*}
\begin{array}{ll} 
   \min \  & x^TQ_0x+2c_0^Tx \\ 
s.t. \  & A_i^Tx=d_i, \ \ \;\;\;\; i=1, ..., p \\
   &  x_i\in \{0, 1\} \;\;\;\;\;\;\;\; i\in \mathcal{I}\\
  &   0\leq x\leq e, 
\end{array}
\end{equation*}
where index set $\mathcal{I}\subseteq \{1, ..., n\}$ denotes binary variables. This problem can be formulated as 
\begin{equation*}
\begin{array}{ll} 
  \min \  & x^TQ_0x+2c_0^Tx \\ 
  s.t. \  &  A_i^Tx=d_i, \ \ \;\;\;\; i=1, ..., p \\
   &  x_i(1-x_i)\leq 0\;\;\;\;\;\;\;\; i\in \mathcal{I}\\
  &   0\leq x\leq e.
\end{array}
\end{equation*}
As all conditions of Theorem \eqref{T1} holds for the above problem, bound \eqref{L4} is exact for  linearly constrained quadratic programs with binary variables.   
\end{proof} 
\end{proposition}
By Theorem 2.6 in \cite{Bur} and strong duality for conic programs $q\in \mathcal{A}_+(B)$ if and only if the following system has a solution 
\begin{align*}
& q(x)+\sum_{i=1}^{n} \alpha_i(x, s)(x_i+s_i-1)\in \mathcal{Q}_+(\mathbb{R}_+^{2n}),\\
& \alpha_i(x, s) \in \mathcal{A}(\mathbb{R}_+^{2n}), \ \  \  i=1, ..., n
\end{align*}
where slack variable $s\in \mathbb{R}^{n}$. A quadratic function $q(x)=x^TQx+2c^Tx+b\in \mathcal{Q}_+(\mathbb{R}_+^{n})$ if and only if matrix $\mathcal{M}(q)$ is copositive \cite{Bom}. So bound \eqref{L4} can be casted as a copostive program. Copostive programs are intractable in general. In fact, they are NP-hard. Nonetheles, there are efficient methods which approximate copositive cone \cite{Bomm, Parrilo}. \\
\indent
It is well-known Shor relaxation is the dual of \eqref{P} when affine multipliers are constant functions. In addition, we have shown before, the dual of \eqref{L2} is Shor relaxation with partial first-level RLT. It is may be of interest to know whether other SDRs can be also obtained in this manner. In the sequel, we will show that some SDRs can be obtained as the  dual of bounds in the form of \eqref{L2} by a suitable choice of affine multipliers or adding some valid cuts.\\
\indent 
Let $m_c=\{i: Q_i\succeq 0, Q_i\neq 0\}$.  In the rest of the section, we make the assumption that for each $i\in m_c $ there exists $\bar x^i \in\mathbb{R}^n$ such that $(\bar x^i)^TQ_i\bar x^i+2c^T\bar x^i< b_i$. Due to the semi-positiveness of $Q_i$ ($i\in m_c$), there exists matrix $R_i$ with $Q_i=R_i^TR_i$.  By Schur Complement Lemma, 
$x^TQ_ix+2c_i^Tx\leq b_i$ is equivalent to
$$
\begin{pmatrix}
I  &  -R_ix\\  -x^TR_i^T  & -2c_i^Tx+b_i,
\end{pmatrix} \succeq 0. 
$$ 
So affine function $\alpha(x)=f^Tx+g$ is non-negative on $L_i=\{x\in \mathbb{R}^n: x^TQ_ix+2c_i^Tx\leq b_i\}$ if and only if the optimal value of the following semi-definite program is greater  than or equal to $-g$.
\begin{align}\label{SDF}
\nonumber \ & \min \ f^Tx \\
  \   &  s.t.\  
\begin{pmatrix}
I  &  -R_ix\\  -x^TR_i^T  & -2c_i^Tx+b_i,
\end{pmatrix} \succeq 0. 
\end{align}
We remark that problem \eqref{SDF} can be reformulated as a second-order cone program. Recently, Zheng et al. proposed some SDRs for QCQPs \cite{Zheng}. In fact, they introduced a unified framework for generating convex relaxations for  QCQPs. They propose the following SDR for \eqref{P} 
 \begin{align}\label{S3}
\nonumber  \min \ & Q_0\bullet X+2c_0^Tx \\
\nonumber  \    s.t.\  &  Q_i\bullet X+2c_i^Tx\leq b_i, \ \ \;\;\;\;\;\;\;\;\;\;\;\;\;i=1, ..., m \\
\nonumber  \   &  XA_i^T= d_ix, \ \ \;\;\;\;\;\;\;\;\;\;\;\;\;i=1, ..., p \\
\nonumber  \  & Ax=d,  \\
\nonumber &  X\geq 0, \ X-xx^T\succeq 0,\\
  &   ex^T-X\geq 0,\\
\nonumber  &    X-ex^T-xe^T+ee^T\geq 0,\\
\nonumber  &
\begin{pmatrix}
x_kI  &  R_iXe_k\\  (R_iXe_k)^T  &  -2c_i^TXe_k+b_ie_k^Tx
\end{pmatrix}
\succeq 0,  \ \ \;\;\;\;\; i=1, ...,m_c, k=1, ..., n  
\\
\nonumber  &  
\begin{pmatrix}
(1-x_k)I  &  -R_iXe_k+R_ix\\  (-R_iXe_k+R_ix)^T  &  2c_i^TXe_k-(2c_i^T+b_ie_k^T)x+b_k
\end{pmatrix}
\succeq 0,    i=1, ...,m_c, k=1, ..., n  
\end{align}
and they called it SDP relaxation with rank-2 second-order cone valid inequalities. Note that the above SDR is obtained by adding  the last two constraints of \eqref{S3} to Shor relaxation with partial first-level RLT. We demonstrate that SDR \eqref{S3} is the dual of  the following bound 
\begin{align}\label{B2}
\nonumber \ & \max \ \ell \\
\nonumber  \   &  s.t.\   x^TQ_0x+2c^T_0x-\ell+\sum_{i=1}^{m} \lambda_i(x^TQ_ix+2c^T_ix-b_i)+ \sum_{i=1}^{p} \alpha_i(x)(A_ix-d_i)-\\
\nonumber  & \; \; \; \; \; \; \; \; \; \; \; \;         \sum_{i=1}^{n} (\beta_i(x)+ \sum_{j\in m_c}\beta_{ij}(x))x_i+\sum_{i=1}^{n} (\gamma_i(x)+ \sum_{j\in m_c}\gamma_{ij}(x))(x_i-1)\in \mathcal{Q}_+(\mathbb{R}^{n}),\\
& \; \; \; \; \; \; \; \; \; \; \; \;     \ \alpha_i\in \mathcal{A}(\mathbb{R}^n), \ i=1, ..., p\\
\nonumber & \; \; \; \; \; \; \; \; \; \; \; \;   \lambda\geq 0,\ \beta_i, \gamma_i \in \mathcal{A}_+(\mathcal{X}), \ i=1, ..., n\\
\nonumber & \; \; \; \; \; \; \; \; \; \; \; \;   \beta_{ij}, \gamma_{ij} \in \mathcal{A}_+(L_j), \ i=1, ..., n, j\in m_c
\end{align}
It is easily seen \eqref{B2} is a bound for \eqref{P}.  For convenience, to show  bound  \eqref{B2} is the dual of \eqref{S3} we consider the  QCQP 
\begin{equation}\label{PS}
\begin{array}{ll} 
 \min \  & x^TQ_0x+2c_0^Tx \\ 
   s.t. \  & x^TQ_1x+2c_1^Tx\leq b_1, \\
  & a_1^Tx\leq d_1,
  \end{array} 
\end{equation}
which has a convex quadratic constraint and a linear inequality constraint. Bound  \eqref{B2} for problem \eqref{PS} is formulated as  
\begin{equation}\label{Exrr}
\begin{array}{ll} 
 \ & \max \  \ell \\ 
 &  s.t. \  x^TQ_0x+2c_0^Tx-\ell +\alpha (x^TQ_1x+2c_1^Tx- b_1)+( f^Tx+g)( a_1^Tx- d_1)\in  \mathcal{Q}_+(\mathbb{R}^{n})\\
  & \; \; \; \; \; \; \; \; \; \alpha\geq 0, \  f^Tx+g\in  \mathcal{A}_+(L_1)
  \end{array} 
\end{equation}
As $\Int(L_1)\neq \emptyset$, the strong duality holds for problem \eqref{SDF}. Accordingly, $f^Tx+g\in  \mathcal{A}_+(L_1)$ is equivalent that the optimal value of the following semi-definite program is   greater or equal to $-g$,
\begin{align*}
  \max \ & -I \bullet Y-b_1y_0\\
   s.t.\  & -R_1y-c_1y_0=\frac{1}{2} f \\
  &
\begin{pmatrix}
Y &  y\\  y^T  & y_0
\end{pmatrix} \succeq 0. 
\end{align*}
Hence,  problem \eqref{Exrr} is reformulated as 
 \begin{align*}
  \max \ &  \ell \\ 
 s.t. \   &  x^TQ_0x+2c_0^Tx-\ell +\alpha (x^TQ_1x+2c_1^Tx- b_1)+\\
 & \; \; \; \;  (2(-R_1y-c_1y_0)^Tx+g)( a_1^Tx- d_1)\in  \mathcal{Q}_+(\mathbb{R}^{n})\\ 
  & -I \bullet Y-b_1y_0\geq -g\\
   & \alpha\geq 0, \  \begin{pmatrix}
Y &  y\\  y^T  & y_0
\end{pmatrix} \succeq 0
\end{align*}
By a little algebra, the dual of the above problem is formulated as 
 \begin{align*}
  \min \   & Q_0\bullet X+2c_0^Tx \\ 
   s.t. \   & Q_1\bullet X+2c_1^Tx\leq b_1,\\
 &    a_1^Tx\leq d_1,\\
 & X-xx^T\succeq 0\\
 & \begin{pmatrix}
(d_1-a_1^Tx)I  &  -R_1Xa_1+d_1R_1x\\  (-R_1Xa_1+d_1R_1x)^T  &  2c_1^TXa_1-(2d_1c_1+b_1a_1)^Tx+b_1d_1
\end{pmatrix}
 \succeq 0,
\end{align*}
which clarifies the point that  bound \eqref{B2} is the dual of  \eqref{S3}. Since \eqref{S3} is strongly feasible, strong duality also holds.\\
\indent
A typical method to tighten the relaxation gap is adding valid cuts. Zheng et al. introduced a class of quadratic valid cuts for QCQP and  they proposed a new SDR by using these valid cuts \cite{Zheng}. Their method generate a quadratic valid cut as follows. Let $F\subseteq \Omega$. Suppose that $u\in \Int(\mathbb{R}^n_+)$ and $0<u_\Omega=\{\max \ u^Tx: \ x\in\Omega\}$. They showed that for $S\succeq 0$, the convex quadratic inequality $x^TSx-u_\Omega \Diag(S)^T\diag(u)^{-1}x\leq 0$ is valid for \eqref{P}, see Proposition 3 in \cite{Zheng}. We remark that the set of generated cuts by this method forms a convex cone in $\mathcal{Q}({\mathbb{R}^n})$.\\
\indent
By the above discussion, one can extend bound \eqref{B2} as follows,
\begin{align*}
\nonumber \ & \max \ \ell \\
\nonumber  \   &  s.t.\   x^TQ_0x+2c^T_0x-\ell+\sum_{i=1}^{m} \lambda_i(x^TQ_ix+2c^T_ix-b_i)+ \sum_{i=1}^{p} \alpha_i(x)(A_ix-d_i)-\\
\nonumber  & \; \; \; \; \; \; \; \; \; \; \; \;         \sum_{i=1}^{n} (\beta_i(x)+ \sum_{j\in m_c}\beta_{ij}(x))x_i+\sum_{i=1}^{n} (\gamma_i(x)+ \sum_{j\in m_c}\gamma_{ij}(x))(x_i-1)+\\
& \; \; \; \; \; \; \; \; \; \; \; \; \sum_{i\in \mathbb{R}} \mu_i(x^TS_ix-u_\Omega \Diag(S_i)^T\diag(u)^{-1}x) \in \mathcal{Q}_+(\mathbb{R}^{n}),\\
& \; \; \; \; \; \; \; \; \; \; \; \;     \ \alpha_i\in \mathcal{A}(\mathbb{R}^n), \ i=1, ..., p\\
\nonumber & \; \; \; \; \; \; \; \; \; \; \; \;   \lambda\geq 0,\ \beta_i, \gamma_i \in \mathcal{A}_+(\mathcal{X}), \ i=1, ..., n\\
\nonumber & \; \; \; \; \; \; \; \; \; \; \; \;   \beta_{ij}, \gamma_{ij} \in \mathcal{A}_+(L_j), \ i=1, ..., n, j\in m_c\\
\nonumber & \; \; \; \; \; \; \; \; \; \; \; \;   \mu_i\geq 0,\ S_i \succeq 0, \  i\in\mathbb{R} 
\end{align*}
which is a non-convex optimization problem with infinite constraints and variables. As mentioned above, the set of valid cuts is a convex cone, so the above bound can be formulated  as semi-definite program
\begin{align}\label{B4}
\nonumber \  \max \ & \ell \\
\nonumber  \     s.t.\  & x^TQ_0x+2c^T_0x-\ell+\sum_{i=1}^{m} \lambda_i(x^TQ_ix+2c^T_ix-b_i)+ \sum_{i=1}^{p} \alpha_i(x)(A_ix-d_i)-\\
\nonumber  &  \; \; \;  \; \;           \sum_{i=1}^{n} (\beta_i(x)+ \sum_{j\in m_c}\beta_{ij}(x))x_i+\sum_{i=1}^{n} (\gamma_i(x)+ \sum_{j\in m_c}\gamma_{ij}(x))(x_i-1)+\\
& \nonumber   \; \; \;   \; \;  x^TSx-u_\Omega \Diag(S)^T\diag(u)^{-1}x \in \mathcal{Q}_+(\mathbb{R}^{n}),\\
& \alpha_i\in \mathcal{A}(\mathbb{R}^n), \ i=1, ..., p\\
\nonumber &  \lambda\geq 0,\ S \succeq 0, \beta_i, \gamma_i \in \mathcal{A}_+(\mathcal{X}), \ i=1, ..., n\\
\nonumber &   \beta_{ij}, \gamma_{ij} \in \mathcal{A}_+(L_j), \ i=1, ..., n, j\in m_c
\end{align}
The dual of bound \eqref{B4} is formulated as 
 \begin{align}\label{S4}
\nonumber  \min \ & Q_0\bullet X+2c_0^Tx \\
\nonumber  \    s.t.\  &  Q_i\bullet X+2c_i^Tx\leq b_i, \ \ \;\;\;\;\;\;\;\;\;\;\;\;\;i=1, ..., m \\
\nonumber  \   &  XA_i^T= d_ix, \ \ \;\;\;\;\;\;\;\;\;\;\;\;\;i=1, ..., p \\
\nonumber  \  & Ax=d,  \\
\nonumber &  X\geq 0, \ X-xx^T\succeq 0,\\
  &   ex^T-X\geq 0,\\
\nonumber  &    X-ex^T-xe^T+ee^T\geq 0,\\
\nonumber  &
\begin{pmatrix}
x_kI  &  R_iXe_k\\  (R_iXe_k)^T  &  -2c_i^TXe_k+b_ie_k^Tx
\end{pmatrix}
\succeq 0,  \ \ \;\;\;\;\; i=1, ...,m_c, k=1, ..., n  
\\
\nonumber  &  
\begin{pmatrix}
(1-x_k)I  &  -R_iXe_k+R_ix\\  (-R_iXe_k+R_ix)^T  &  2c_i^TXe_k-(2c_i^T+b_ie_k^T)x+b_k
\end{pmatrix}
\succeq 0,   \ \  i=1, ...,m_c, k=1, ..., n  \\
\nonumber &  u_\Omega \diag(u)^{-1}\diag(x)-X \succeq 0,
\end{align}
which is the SDR proposed in \cite{Zheng}; See $SDP_{\alpha_u}$. Since \eqref{S4} is also strongly feasible, we have strong duality. Here we just investigate some well-known SDRs and show that they can be interpreted as the dual of a bound in the form \eqref{L2}. However, by a similar argument one can show that most SDRs can be obtained as the dual of a bound in the form of \eqref{L2}.\\
\indent
We conclude the section by mentioning some points. As the dual of the proposed bounds are well-known SDRs, we have just reinvented the wheel. Of course, this statement is correct, but viewing SDRs from this aspect can supply us with more tools for analyzing a SDR method. For instance, one can extend bound \eqref{B2} as follows
\begin{align}\label{B5}
\nonumber \  \max \ & \ell \\
\nonumber  \     s.t.\ &  x^TQ_0x+2c^T_0x-\ell+\sum_{i=1}^{k} \lambda_i(x^TQ_ix+2c^T_ix-b_i)+ \sum_{i=1}^{k} \alpha_i(x)(A_ix-d_i)-\\
  & \; \; \; \; \;       \sum_{i=1}^{k} \beta_i(x)x_i+\sum_{i=1}^{k} \gamma_i(x)(x_i-1)\in \mathcal{Q}_+(\mathbb{R}^{n}),\\
\nonumber &     \ \alpha_i\in \mathcal{A}(\mathbb{R}^n), \ i=1, ..., p\\
\nonumber &    \lambda_i \geq 0,\ \beta_i, \gamma_i \in \mathcal{A}_+(V), \ i=1, ..., n
\end{align}
where $V=\{x\in \mathcal{X}: x^TQ_ix+2c_i^Tx\leq b_i, \ i\in m_c\}$. As $ \mathcal{A}_+(\mathcal{X})\cup_{i\in m_c} \mathcal{A}_+(L_i)\subseteq  \mathcal{A}_+(V)$, bound \eqref{B5} dominates   \eqref{B2}. Therefore, the dual of   \eqref{B5} leads to a SDR which dominates  \eqref{S3}. Here, it is assumed that $\Int(\cap_{i\in m_c} L_i)\cap \mathcal{X}\neq \emptyset$. In the same line, we can formulate the following bound which dominates  \eqref{B4}
\begin{align}\label{B6}
\nonumber \  \max \ & \ell \\
\nonumber  \     s.t.\  & x^TQ_0x+2c^T_0x-\ell+\sum_{i=1}^{k} \lambda_i(x^TQ_ix+2c^T_ix-b_i)+ \sum_{i=1}^{k} \alpha_i(x)(A_ix-d_i)-\\
\nonumber   & \; \; \; \; \;       \sum_{i=1}^{k} \beta_i(x)x_i+\sum_{i=1}^{k} \gamma_i(x)(x_i-1)+ x^TSx-u_\Omega \Diag(S)^T\diag(u)^{-1}x \in \mathcal{Q}_+(\mathbb{R}^{n}),\\
 & \alpha_i\in \mathcal{A}(\mathbb{R}^n), \ i=1, ..., p\\
\nonumber &  \lambda\geq 0,\ S \succeq 0, \beta_i, \gamma_i \in \mathcal{A}_+(W), \ i=1, ..., n
\end{align}
where $W=\{x\in V: x^TSx-u_\Omega \Diag(S)^T\diag(u)^{-1}x\leq 0,\  \forall S\succeq 0\}$. Note that if there exists $\bar x\in \mathcal{X}$ with $\bar x^TS\bar x-u_\Omega \Diag(S)^T\diag(u)^{-1}\bar x < 0$ for $S\succeq 0$ and $\bar x^TQ_i\bar x+2c_i^T\bar x< b_i$ ($ i\in m_c$), then $f^Tx+g\in\mathcal{A}_+(W)$ is equivalent to the consistency of the system 
\begin{align*}
& f^Tx+g+\lambda^T(Ax-d)+\mu^T(x-e)-\nu^Tx+x^TSx-u_\Omega \Diag(S)^T\diag(u)^{-1}x\in \mathcal{Q}_+(\mathbb{R}^{n}), \\
& \mu, \nu\geq 0, \ S\succeq 0.
\end{align*}
Thus, bound \eqref{B6} is reformulated as a semi-definite program, and consequently its dual gives a SDR which dominates \eqref{S4}.\\
\indent
Another point about the proposed bounds is that they not only provide a lower bound, but also give a convex underestimator. The given convex underestimator can be employed in optimization methods for generating a solution. \\
\indent
 It is well-known when an optimal solution of a SDR has rank one the SDR is exact \cite{Luo}. The next proposition gives necessary and sufficient conditions for exactness.  For convenience to state the proposition, we consider bound \eqref{L2}. Let 
$F_{opt}$ denote the optimal solution set of  \eqref{P}.
 \begin{proposition}\label{p5}
Bound  \eqref{L2} is exact if and only if there exists feasible point $\bar \lambda$, $\bar \alpha_i$ ($i=1, ..., p$), $\bar \beta_i, \bar \gamma_i$ ($i=1, ..., n$) and $\bar \ell$ 
  with
\begin{align*}
\bar x\in \argmin_{x\in \mathbb{R}^n} & \ x^TQ_0x+2c^T_0x-\bar\ell+\sum_{i=1}^{m} \bar\lambda_i(x^TQ_ix+2c^T_ix-b_i)+ \\
 & \; \; \; \;    \sum_{i=1}^{p} \bar\alpha_i(x)(A_ix-d_i)-\sum_{i=1}^{n} \bar\beta_i(x)x_i+\sum_{i=1}^{n} \bar\gamma_i(x)(x_i-1)\\
 & \bar x^TQ_0\bar x+2c^T_0\bar x=\bar \ell, \  \forall\bar x \in F_{opt}
 \end{align*}
\begin{proof} 
Let \eqref{L2} be exact and suppose that $\bar \lambda$, $\bar \alpha_i$ ($i=1, ..., p$), $\bar \beta_i, \bar \gamma_i$ ($i=1, ..., n$) and $\bar \ell$ is an optimal solution. As the bound is exact, we have $\bar x^TQ_0\bar x+2c^T_0\bar x=\bar \ell$ for  $\bar x\in F_{opt}$. In the light of $q(x)=x^TQ_0x+2c^T_0x-\bar\ell+\sum_{i=1}^{k} \bar\lambda_i(x^TQ_ix+2c^T_ix-b_i)+\sum_{i=1}^{k} \bar\alpha_i(x)(A_ix-d_i)-\sum_{i=1}^{k} \bar\beta_i(x)x_i+\sum_{i=1}^{k} \bar\gamma_i(x)(x_i-1)\in\mathcal{Q}_+(\mathbb{R}^{n})$ and $q(\bar x)=0$, we have
$$
\bar x \in\argmin\{q(x): x\in \mathbb{R}^{n}\}, 
$$
which completes the if part. The only-if part is immediate. 
\end{proof} 
\end{proposition}
It is worth mentioning that as strong duality holds for all proposed bounds,  exactness of SDRs and bounds are equivalent. Moreover, bound \eqref{L2} or SDR \eqref{D2} are exact for general QCQP if and only if $n=2$, see \cite{anst} for more details.
\section{Cubic bounds}
\noindent
In this section, we propose cubic bounds for QCQP. So far, we use affine functions as dual variables. The most important point for applying other functions is that the obtained problem should be tractable.\\
\indent
 Due to the structure of  \eqref{P}, one may consider the following convex cones for dual variables,
 \begin{enumerate}
 \item 
 $\mathcal{Q}_+^c(\mathcal{X})$: non-negative convex quadratic functions on $\mathcal{X}$,
  \item 
   $\mathcal{Q}^N(\mathbb{R}^{n})$: quadratic functions with non-negative coefficients.
  \end{enumerate}
  Both above-mentioned cones have non-empty interior and verifying the membership of a given quadratic function is tractable. Verifying $q\in \mathcal{C}^N(\mathbb{R}^n)$ is straightforward. By alternative theorem \cite{Mang}, $q(x)=x^T\hat Q x+2\hat c^T c+\hat c_0$ belongs to $\mathcal{Q}_+^c(\mathcal{X})$ if and only if there exist $\lambda\in \mathbb{R}^{p}$ and 
  $\mu, \nu_2\in \mathbb{R}^{n}_+$ with 
  $$
 \mathcal{M}(x^T\hat Q x+2\hat c^T c+\hat c_0+\lambda^T(Ax-d)+\mu^T(x-e)-\nu^Tx)\succeq 0.
  $$
  
  By employing quadratic functions as dual variables, we are faced with checking non-negativity of a cubic function. Of course, a cubic function may not be non-negative on $\mathbb{R}^n$, unless it is quadratic. So it appears our effort by substituting some class of quadratic functions instead of affine functions was in vain.
Nevertheless,  checking non-negativity of some class of cubic functions are tractable on non-negative orthant. For instance, one may consider the following sets of  cubic functions,
 \begin{enumerate}
 \item 
 $\mathcal{C}_+^c(\mathbb{R}_{+}^n)$: non-negative convex cubic functions on $\mathbb{R}_{+}^n$,
  \item 
   $\mathcal{C}^N(\mathbb{R}^n)$: cubic functions with non-negative coefficients,
  \end{enumerate}
  Both sets are convex cones with  non-empty interior. In addition, to check a cubic function belongs to these cones are tractable. Let $q(x)=Tx^3+xQx+cx+c_0$ be a cubic function, where $T$ is a symmetric tensor of order 3. Verifying $q\in \mathcal{C}^N(\mathbb{R}^n)$ is straightforward. To check $q\in\mathcal{C}_+^c(\mathbb{R}_{+}^n)$, we need first to impose the following conditions 
  \begin{equation*}
\begin{array}{ll} 
Te_i\succeq  0, \  i=1, ..., n,\\
Q\succeq  0,
\end{array}
\end{equation*}
which guarantees convexity of $q$ on $\mathbb{R}_{+}^n$. As $q$ is convex, its optimal value can be obtained by primal interior point methods. As a result, membership verification is tractable in this case, but cannot be checked explicitly by some linear inequalities.\\
\indent
Another replacement for $\mathcal{C}_+^c(\mathbb{R}_{+}^n)$ or $\mathcal{C}^N(\mathbb{R}^n)$ may be the set of quadratically Sum-of-Squares. We call a cubic function $q(x)=Tx^3+x^TQx+c^Tx+c_0$ quadratically Sum-of-Squares if  $T(x^{(2)})^3+(x^{(2)})^TQ(x^{(2)})+c(x^{(2)})+c_0$ is Sum-of-Squares, where $x^{(2)}=(x_1^2, ...,  x_n^2)$. Note that one can check whether a polynomial is Sum-of-Squares using semidefinite programming \cite{Las}.\\
\indent
By the above discussion, we propose the following bound for  \eqref{P}
 \begin{align}\label{B7}
\nonumber \  \max \ & \ell \\
\nonumber  \     s.t.\   & x^TQ_0x+2c^T_0x-\ell+\sum_{i=1}^{m} \lambda_i(x)(x^TQ_ix+2c^T_ix-b_i)+\sum_{i=1}^{p} \alpha_i(x)(A_ix-d_i)- \\
\nonumber  & \; \; \; \;         \sum_{i=1}^{n} \beta_i(x)x_i+\sum_{i=1}^{n} \gamma_i(x)(x_i-1)-\kappa(x)\in \mathcal{Q}_+(\mathbb{R}^n),\\
&    \lambda_i\in\mathcal{A}_+(X), \ \alpha_j\in \mathcal{Q}(\mathbb{R}^n), \ i=1, ..., m, \ j=1, ..., p\\
\nonumber & \beta_i, \gamma_i \in \mathcal{Q}_+^c(\mathcal{X}), \ i=1, ..., n\\
\nonumber  & \kappa\in\mathcal{C}^N(\mathbb{R}^n).
\end{align}
By the above discussion, problem \eqref{B7} can be formulated as a semi-definite program and it has $O(n^3)$ variables. By Proposition \ref{p1}, bound \eqref{B7} is always finite and generates a lower bound greater or equal to that of  \eqref{L2}.\\
\indent
One may wonder how the bound given by \eqref{B7} can be improved. The straightforward method for tightening can be enlargement of feasible set. In problem \eqref{B7}, we have  linear and quadratic function variables. One can adopt methods in Section 2 to tighten bound \eqref{B7}.
 The following proposition gives necessary and sufficient conditions under which bound \eqref{B7} is exact.
 \begin{proposition}\label{pc1}
Bound  \eqref{B7} is exact if and only if there exists feasible point $\bar \lambda_i$ ($i=1, ..., m$), $\bar \alpha_i$ ($i=1, ..., p$), $\bar \beta_i, \bar \gamma_i$ ($i=1, ..., n$), $\bar \ell$ and $\bar\kappa$ 
  with
 \begin{align*}
\bar x\in \argmin_{x\in \mathbb{R}^n} & \ x^TQ_0x+2c^T_0x-\bar\ell+\sum_{i=1}^{m} \bar\lambda_i(x)(x^TQ_ix+2c^T_ix-b_i)+ \\
 & \; \; \; \;    \sum_{i=1}^{p} \bar\alpha_i(x)(A_ix-d_i)-\sum_{i=1}^{n} \bar\beta_i(x)x_i+\sum_{i=1}^{n} \bar\gamma_i(x)(x_i-1)-\bar\kappa(x)\\
 & \bar x^TQ_0\bar x+2c^T_0\bar x=\bar \ell, \ \forall\bar x \in F_{opt}
 \end{align*}
\begin{proof} 
Similar to Proposition \ref{p5} is proved. 
\end{proof} 
\end{proposition}
 In the same line for problem \eqref{B7}, one could consider quadratics or linear multipliers for which optimal value of  \eqref{B7} are non-negative. In fact, one may consider 
 $q(x)=x^T\hat Qx+2\hat c^Tx+\hat c_0$ eligible if  the following system has a solution
  \begin{align*}
\nonumber  \    \   & x^T\hat Qx+2\hat c^Tx+\hat c_0+\sum_{i=1}^{m} \lambda_i(x)(x^TQ_ix+2c^T_ix-b_i)+\sum_{i=1}^{p} \alpha_i(x)(A_ix-d_i)- \\
\nonumber  & \; \; \; \;         \sum_{i=1}^{n} \beta_i(x)x_i+\sum_{i=1}^{n} \gamma_i(x)(x_i-1)-\kappa(x)\in \mathcal{Q}_+(\mathbb{R}^n),\\
&    \lambda_i\in\mathcal{A}_+(X), \ \alpha_j\in \mathcal{Q}(\mathbb{R}^n), \ i=1, ..., m, \ j=1, ..., p\\
\nonumber & \beta_i, \gamma_i \in \mathcal{Q}_+^c(\mathcal{X}), \ i=1, ..., n\\
\nonumber  & \kappa\in\mathcal{C}^N(\mathbb{R}^n).
\end{align*}
 Thus, we obtain a hierarchy for tackling  \eqref{P}. This  hierarchy is increasing and  each problem is formulated as a semi-definite program. If we consider problem \eqref{B7} as a first problem of hierarchy, the $k^{th}$ problem contains  $O(n^{k+2})$ variables. The most important inquiry concerning this method is its convergence in finite steps. In addition, if it is  convergent in finite steps, what the order of $k$ for which we have convergence. As the subject of the paper is quadratic and cubic bounds, we leave these questions for further research.\\
 \indent
 The following example demonstrate that \eqref{B7} could generate a bound tighter than the proposed bounds in Section 2. 
 \begin{example}
Consider the nonconvex QCQP,
 \begin{align*}
\min \ & -8x_1^2-x_2^2+x_3^2-5x_4^2 +14x_1x_2+10x_1x_4+4x_2x_4-10x_2 \\
s.t.\ & 2x_1^2 +2x_2^2+4x_1x_2 +8x_1+6x_2+x_4\leq 8\\
 & -8x_1^2-5x_2^2+2x_1x_4-5x_2^2 -4x_1+4x_2+2x_4 \leq -4,\\
  & 2x_1^2 +x_2^2+x_4^2+2x_1+x_4 \leq 4,\\
& x_1 + 2x_2+2x_3+x_4 = 3, \\
& x\in B.
 \end{align*}
  The problem has two convex quadratic constraints, one nonconvex quadratic constraint and nine linear constraints, with the optimal value is  $-8.0008$ and the optimal solution 
  $ \begin{pmatrix} 0.4203,  &0.4942,  & 0.7956, & 0\end{pmatrix}$. We set $u=\begin{pmatrix} 1, & 2, & 2, & 1 \end{pmatrix}^T$ and $u_\Omega = max\{u^T x: \ x\in F\} = 3.9145$.  The performance of the bounds are listed in Table 1, which lb denotes the generated lower bound.

\begin{table}[htbp]
\centering
\begin{tabular}
{  c  c   c   c   c   c   c  c  }
  \toprule
\shortstack{Bound \\ {\color{white} hted}}       &  \shortstack{Shor \\ relaxation}  & \shortstack{Bound \\ \eqref{L2} }  &    \shortstack{Bound \\ \eqref{B2} } &    \shortstack{Bound \\ \eqref{B4} } &  \shortstack{Bound \\ \eqref{B5} } &   \shortstack{Bound \\ \eqref{B6} } &   \shortstack{Bound \\ \eqref{B7} } \\    
   \hline
  \\[-1ex]   
\multicolumn{1}{l}{lb} &    -44.0945    & -15.2676            & -13.3647             & -13.2518              &   -13.2294         &   -11.8969           &   -8.0008  \\ 
  \hline
    \hline
\end{tabular}
\end{table}
  As seen bound \eqref{B7} is exact for the example.
\end{example}
  In the rest of the section, we investigate the relationship between bound  \eqref{B7} and the conventional lower bounds for QCQPs. Due to the computational burdensome, cubic bounds are not commonplace and they have been applied just for some types of QCQPs such as  standard quadratic programs. Of course, we can provide a comparison between bound \eqref{B7} and general polynomial optimization methods, including Lasserre hierarchy \cite{Las}, with  $O(n^3)$ variables, but we prefer bounds tailored for QCQPs.\\
  \indent
  Consider the standard quadratic program,
\begin{equation}\tag{StQP}\label{SQ}
\begin{array}{ll} 
  \min \ & x^TQx\\
  s.t. \ & \sum _{i=1}^n x_i= 1,\\
   &  x\geq 0.
\end{array}
\end{equation}
It is well-known that \eqref{SQ} is solvable in polynomial time provided $Q$ is either positive semi-definite or negative semi-definite on standard simplex. In general, however, \eqref{SQ} is NP-hard \cite{Bomm}. Suppose that $\Delta$ denote the standard simplex.\\
\indent
Let $\ell_Q$ denote the optimal value of \eqref{SQ}.  
We remark that optimizing a quadratic function on standard simplex can be formulated as \eqref{SQ}. This is resulted from the fact that for each $x\in \Delta$, we have 
$x^TQx+2c^Tx=x^T(Q+ec^T+ce^T)x$.\\
\indent
One effective method for handling (StQP) is Parrilo hierarchy \cite{Bomm}. In this method, for $r = 0, 1, ...$ the following problem gives a lower bound
\begin{equation}\label{Par}
p^r_Q =max\{\ell: Q-\ell ee^T\in \mathcal{P}^r\},
\end{equation}
where $\mathcal{P}^r = \{A : (􏰄\sum_{i=1}^n  x_i^2)(􏰄􏰄\sum_{i=1}^n 􏰄\sum_{j=1}^n  x_i^2A_{ij} x_j^2 )\in \Sigma[x]\}$ and $\Sigma[x]$ denotes set of all sum of square polynomials. It is well-known for sufficiently large $r$, $p^r_Q$ is equal to the optimal value of (StQP) \cite{Bomm}. In addition, the number of variables of \eqref{Par} is of 
$O(n^{r+2})$ \cite{Parrilo}.\\
\indent
Bound \eqref{B7} is formulated as follows for \eqref{SQ},
 \begin{align}\label{BSQ}
\nonumber \  \max \ & \ell \\
\nonumber  \     s.t.\   & x^TQx-\ell -\sum_{i=1}^{n} \alpha_i(x)x_i+ \alpha_{n+1}(x)(e^Tx-1)-\kappa(x)\in \mathcal{Q}_+(\mathbb{R}^n),\\
 & \alpha \in \mathcal{Q}_+^c(\Delta), \ i=1, ..., n\\
\nonumber  & \kappa\in\mathcal{C}^N(\mathbb{R}^n).
\end{align}
As mentioned earlier, the number of variables of \eqref{BSQ} is of $O(n^3)$. So, one may wonder what is the relationship between $p^1_Q$ and the optimal value of $\eqref{BSQ}$ The following theorem says these bounds are equivalent. Before we get to the proof, let us mention some points. It is shown in \cite{Bomm, Parrilo}, the symmetric
matrix $B\in \mathcal{P}^1$ if and only if there exist symmetric matrices $K^{(1)}, ..., K^{(n)}$ with
\begin{align}
& B-K^{(i)} \succeq  0, \  \ & i=1, ..., n\\ \label{CF}
& K^{(i)}_{ii}=  0, \ \  & i=1, ..., n\\
& K^{(j)}_{ii}+2K^{(i)}_{ij}= 0, \  \ & i\neq j\\
& K^{(i)}_{jk}+K^{(j)}_{ik}+K^{(k)}_{ij}\geq 0, \  \ & i> j>k. \label{CE}
\end{align}
Let convex quadratic function $q(x) = x^TSx + 2c^T x + c_0$ be nonnegative on $\Delta$.
It is easily seen that $ (􏰄\sum_{i=1}^n  x_i)^2q(􏰄\sum_{i=1}^n  x_i)^{-1}x)$ is homogeneous polynomial of degree two. So for some symmetric matrix $Q$, we have 
$ (􏰄\sum_{i=1}^n  x_i)^2q(􏰄\sum_{i=1}^n  x_i)^{-1}x)= x^T Qx$. As $q\in\mathcal{Q}_+^c(\Delta)$, there exist nonnegative multipliers $\lambda_i$ ($i=1, ..., n$) and
$\lambda_{n+1}$ with
$$
q(x)-􏰄\sum_{i=1}^n \lambda_i x_i+\lambda_{n+1}(e^Tx-1)\in\mathcal{Q}_+(\mathbb{R}^n).
$$
By the replacement of $x$ with $(\sum_{i=1}^n  x_i)^{-1}x$ and multiplication of $ (􏰄\sum_{i=1}^n  x_i)^2$, it is readily seen that $Q$ can be represented as a summation of a positive semi-definite matrix and a nonnegative matrix. 
 \begin{theorem}
Bounds  $p^1_Q$ and \eqref{BSQ} are equivalent.
\begin{proof} 
 First, we show that the optimal value of $\eqref{BSQ}$ is less than or equal to $p^1_Q$. Let $\bar\alpha_i(x)=x^TS_ix+2c^T_ix+g_i$, $i=1,...,n+1$, $\bar\kappa$ and $\bar\ell$ be an optimal solution of $\eqref{BSQ}$. (Without loss of generality, it is assumed $\eqref{BSQ}$ attains its optimal solution.)  As $\bar\kappa$ is nonnegative on standard simplex,  
 $(e^T x)^3\bar\kappa((e^T x)^{-1}x)$ is homogeneous polynomial of degree three with nonnegative coefficients. Thus, for nonnegative symmetric matrix $K_i$, $i = 1, ..., n$, we have $(e^T x)^3\bar\kappa((e^T x)^{-1}x)=\sum_{i=1}^n x_i(x^T K_ix)$. Furthermore, $(e^T x)^2\bar\alpha_i((e^T x)^{-1}x)=x^T (L_i + M_i)x$, $i = 1, ..., n$, where $L_i$ and $M_i$ are non-negative and positive semi-definite, respectively. Therefore, we have
 $$
 \sum_{i=1}^n x_ix^T(Q-\bar\ell ee^T-K_i-L_i-M_i-L_0-M_0)x=0,
$$
where $L_0\geq 0$ and $M_0\succeq 0$.
As nonnegative diagonal matrices are positive semi-definite and convexity of nonnegative and positive semi-definite matrices, with a little algebra, we get symmetric matrices $\bar K_i$ , $i=1,...,n$, which satisfy \eqref{CF}-\eqref{CE}. As a result, $\bar\ell\leq p^1_Q$.\\
 Now, we prove that $p^1_Q$ is less than or equal to the optimal value of  \eqref{BSQ}. Similar to the former case, we assume that optimal solution is attained. So,  there exist symmetric matrices $K_i$, $i = 1,...,n$, which satisfy  \eqref{CF}-\eqref{CE} for $B = Q-p^1_Q ee^T$. Let $M_i = Q-p^1_Q ee^T-K_i$, $i = 1,...,n$. We have
 $$
 (e^Tx)(x^T(Q-p^1_Q ee^T)x)-\sum_{i=1}^n x_ix^TM_ix\in\mathcal{C}^N(\mathbb{R}^n).
 $$
Hence,
 $$
x^TQx-p^1_Q-\sum_{i=1}^n x_ix^TM_ix+(e^Tx-1)(x^T(Q-p^1_Q ee^T)x)\in\mathcal{C}^N(\mathbb{R}^n).
 $$
which completes the proof.
\end{proof} 
\end{theorem}
The proof of above theorem reveals for bound  \eqref{BSQ}, we can replace  $\mathcal{Q}_+^c(\Delta)$ by the set of homogeneous convex quadratics. So, \eqref{BSQ} is equivalent to the problem
 \begin{align*}
\nonumber \  \max \ & \ell \\
\nonumber  \     s.t.\   & x^TQx-\ell -\sum_{i=1}^{n} (x^TS_ix)x_i+ \alpha_{n+1}(x)(e^Tx-1)-\kappa(x)\in \mathcal{Q}_+(\mathbb{R}^n),\\
 & S_i\succeq 0, \ i=1, ..., n\\
\nonumber  & \kappa\in\mathcal{C}^N(\mathbb{R}^n).
\end{align*}

We conclude the section by noting that bound \eqref{B7} dominates semidefinite relaxations obtained in the Laserre Hierarchy with the same order of variables for QCQPs. Strictly speaking, as  bound \eqref{B7} applies $\mathcal{Q}_+^c(\mathcal{X})$ instead of $\mathcal{Q}_+^c(\mathbb{R}^n)$, it can  generate tighter bounds in comparison with the Laserre Hierarchy with the same order of variables. Moreover,  bound \eqref{B7} and RLT-level 2 are not necessarily relevant. 


 \section{Computational results} 
 \noindent
 As mentioned above, two important factors which determine the efficiency of a given bound are the quality of the generated bound and its computational time. In the section, we compare quadratic bound $\eqref{L2}$ and cubic bound$\eqref{B7}$. The aim of this section is not to compare the numerical performance all above-mentioned bounds comprehensively. We merely present numerical performance a quadratic and a cubic bound on some instances. 
 The reason why we chose these bounds is that they are formulated in the same line with different types multipliers. \\
 \indent
 In order to evaluate the efficiency of these bounds, we generated  26 random QCQPs in the form \eqref{P} in $\mathbb{R}^{20}$ with five non-convex quadratic constraints and two equality constraints. To solve semi-definite programs, we employed MOSEK in Matlab 2018b environment \cite{MOSEK}.  We applied YALMIP to pass bounds  $\eqref{L2}$ and
  $\eqref{B7}$ to MOSEK \cite{yalmip}. In addition, we employed BARON to obtain optimal values \cite{BARON}. Moreover, the computations were run on a Windows PC with Intel Core i7 CPU, 3.4 GHz, and 24GB of RAM. \\
  \indent
  Figure 1 shows the generated bounds via $\eqref{L2}$ and $\eqref{B7}$. In the figure, red line denotes the difference between optimal value and generated bound  by $\eqref{L2}$, absolute gap, and  blue denotes absolute gap for bound $\eqref{B7}$. Bound \eqref{B7} were exact for nine examples. However, on average, the computational time for bounds $\eqref{L2}$ and $\eqref{B7}$ were $0.1$ and $3.8$ seconds, respectively.\\
\begin{center}
\includegraphics[width=0.8\textwidth]{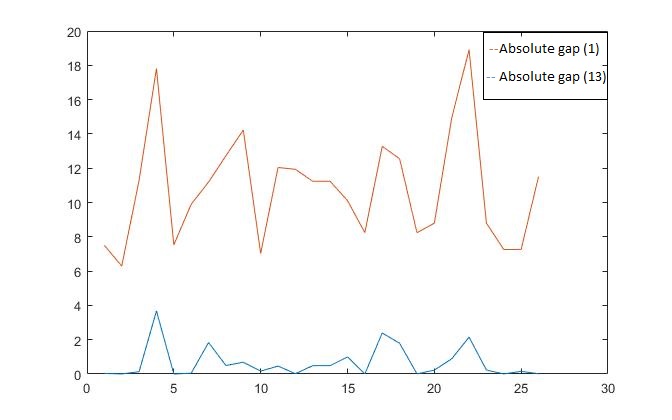}
\captionof{figure}{Absolute gap}
\end{center}

\bibliographystyle{tfs}      
\bibliography{temr}

\end{document}